\documentstyle[amssymb,amsmath,amsthm,draft,epsf,12pt]{article}
\scrollmode
\newcommand{\col}{\mathop{\rm col}\nolimits}
\newcommand{\aC}{{\mathbb C}}
\newcommand{\aR}{{\mathbb R}}
\newcommand{\aN}{{\mathbb N}}

\newcommand{\pd}[2]{\dfrac{\partial #1}{\partial #2}}
\newtheorem{theo}{Theorem}
\newtheorem{lemm}{Lemma}
\newtheorem{cond}{Condition}
\theoremstyle{remark}
\newtheorem{rema}{Remark}
\theoremstyle{definition}

\DeclareMathOperator{\diam}{diam}

\newcommand{\1}{$(A)$}

 \makeatother

\begin{document}
{\large

\begin{center}

{\bf Kryzhevich S.G.}

Saint-Petersburg State University, Russia,

E-mail: kryzhevitz@rambler.ru, kryzhevich@hotmail.com

\bigskip

{\bf Smale Horseshoe and Grazing Bifurcation in Impact Systems}
\medskip

\end{center}}

\noindent\textbf{Abstract.}{ Bifurcations of dynamical systems, described by a second order differential equations and by an impact condition are studied. It is shown that the variation of parameters when the number of impacts of a periodic solution increases, leads to the occurrence of a hyperbolic chaotic invariant set.
}

\textbf{Introduction.} The vibro-impact systems appear in different mechanical problems (modeling of impact dampers, clock mechanisms, immersion of constructions, etc.). All the impact systems are strongly nonlinear. Their properties resemble a lot ones of classical nonlinear systems. Particularly, the chaotic dynamics is possible \cite{ach} -- \cite{31.}, \cite{37.}~-- \cite{45.}, \cite{5.}~-- \cite{62.}, \cite{77.}, \cite{80.}.

There is a big number of publications, devoted to bifurcations, proper to vibro-impact systems. One of them, the so-called grazing bifurcation, first described by A. Nordmark \cite{62.}, corresponds to the case, when there is a family of periodic solution, which has a finite number of impacts over the period and this number increases or decreases as the parameter changes. For the bifurcation value of the parameter, the periodic solution has an impact with a zero normal velocity. It was shown that this bifurcation implies the nonsmooth behavior of solutions, instability of the periodic solution in the parametric neighborhood of grazing and, in additional assumptions, the chaotic dynamics (see \cite{22.}, \cite{31.}, \cite{37.}--\cite{48.}, \cite{62.} and \cite{80.} and the references therein). However, there was no analytical conditions, sufficient for the existence of chaotic invariant sets. For the s.d.f. case these conditions have been obtained in the author's work \cite{6.6.}. In this paper generalize this result for systems with several degrees of freedom. The main result of this paper is the method, which allows to find homoclinic points, corresponding to grazing. The main idea of the proof is the nonsmoothness of Perron surfaces in the neighborhood of periodic solution. If these manifolds bend in a "good"\ way (the corresponding sufficient conditions can be written down explicitly) they can intersect. This implies chaos. We study a motion of a point mass, described by system of second order differential equations of the general form and impact conditions of Newtonian type. In order to avoid technical troubles we assume that the delimiter is plain, immobile and slippery. However, there is no obstacles to apply the offered method to the systems with a mobile delimiter (see, for example \cite{45.}), ones with non-Newtonian model of impacts \cite{121.}, \cite{37.}, \cite{48.}, \cite{4.} and even to some special cases of strongly nonlinear dynamical systems without any impact conditions.

\textbf{1. Mathematical model.} Consider a segment $J=[0,\mu^*]$ and a $C^2$ smooth function $f(t,z,\mu):\aR^{2n+1}\times J\to \aR^n$. Suppose that $f(t,z,\mu)\equiv f(t+T(\mu),z,\mu)$. Here the period $T(\mu)$ is a $C^2$ smooth function of the parameter $\mu$, and $T(0)>0$. We may suppose without loss of generality that $T(\mu)$ does not depend on $\mu$, making, if necessary, the transformation $t'=tT(0)/T(\mu)$. Denote
$$\begin{array}{c}
z=\begin{pmatrix}z_{1} \\ \dots \\ z_{n}\end{pmatrix}, \quad
z_k=\begin{pmatrix}x_k\\ y_k\end{pmatrix},\quad k=1,\ldots, n, \quad
x=\begin{pmatrix}  x_{1} \\ \dots \\ x_{n} \end{pmatrix},\quad \\[5pt]
y=\begin{pmatrix}  y_{1} \\ \dots \\ y_{n} \end{pmatrix},
\quad f=\begin{pmatrix}  f_{1} \\ \dots \\ f_{n}  \end{pmatrix},\quad {\bar x}=\begin{pmatrix}  x_{2} \\ \dots \\ x_{n} \end{pmatrix},\quad
{\bar y}=\begin{pmatrix}  y_{2} \\ \dots \\ y_{n} \end{pmatrix}, \quad
{\bar z}=\begin{pmatrix}  z_{2} \\ \dots \\ z_{n} \end{pmatrix}.\end{array}$$
Consider the system of $2n$ first order ordinary differential equations
\begin{equation}\label{mgr1.1} \dot x_k=y_k;\qquad \dot y_k=f_k(t,z,\mu), \quad k=1,\ldots, n. \end{equation}

Let $r=r(y_1,\mu)\in (0,1]$ be a $C^2$~-- smooth function. Let us note by $\col (a_1,\ldots,a_m)$ the column vector, consisting of elements $a_1,\ldots, a_m$. Suppose that the system \eqref{mgr1.1} is defined for $z\in\Lambda=[0,+\infty)\times {\mathbb R}^{2n-1}$ and the following Newtonian impact conditions take place as soon as the first component $x_1$ of the solution vanishes.

\begin{cond}\label{cog1} \
\begin{enumerate}
\item If $x_1(t_0)=0$ then
$x(t_0+0)=x(t_0-0)$, $y_1(t_0+0)=-r(-y_1(t_0-0),\mu)y_1(t_0-0)$, $\bar y(t_0+0)=\bar y(t_0-0)$.
\item Let the solution $z(t)=\col(z_1(t),\ldots,z_n(t))$ of the system \eqref{mgr1.1} be such that $z_1(t_0)=0$ for a certain instant $t_0$ and $f_1(t_0,0,\bar z(t_0),\mu)\leqslant 0$. Let $\zeta(t)$ be the solution of the system
$$\dot x_k=y_k;\qquad \dot y_k=f_k(t,0,0, \bar z,\mu), \quad k=2,\ldots, n$$
with the initial data $\zeta(t_0)={\bar z}(t_0)$ and the segment $I=[t_0,t_1]$ be such that
$$f_1(t,0,0,\zeta(t),\mu)\leqslant 0$$ for all $t\in I$. Then $z_1(t)=0$ and $\bar z(t)=\zeta(t)$ for any $t\in I$. \end{enumerate}\end{cond}

\begin{rema} Similarly one may consider a domain $\Lambda=\Lambda_x\times\aR^{n}$, where the domain $\Lambda_x$ is bounded by a $C^2$ smooth manifold $M$ of the dimension $n-1$.
\end{rema}

Let us denote the dynamical system, defined by Eq. \eqref{mgr1.1} and mentioned impact conditions by \1.

\noindent\textbf{2. The family of periodic solutions.} Since the solutions of the vibro-impact systems are discontinuous at the impact instants, the classical results on integral continuity are not applicable. Nevertheless, the following two statements hold true.

\begin{lemm}\label{leg1} Let $z(t)=\col (x_1(t),y_1(t),\ldots,x_n(t),y_n(t))$ be the solution of the system \1 for $\mu=\mu_0$ and the initial data $z(t_0)=z^0=\col (x_1^0,y_1^0,\ldots,x_n^0,y_n^0)$ and $x_1^0\neq 0$. Suppose that this solution is defined on the segment $[t_-,t_+]\ni t_0$. Assume that there are exactly $N$ zeros $t_-<\tau^0_1<\ldots<\tau^0_N<t_+$ of the function $x_1(t)$ over the segment $[t_-,t_+]$ and $y_1(\tau^0_j-0)\neq 0$, $(j=1,\ldots,N)$. Then for any $\varepsilon>0$ there exists a neighborhood $U$ of the point $\col (t_0,z^0,\mu_0)$ such that
for any fixed $t\in [t_-,t_+]\setminus\bigcup\limits_{k=1}^N (\tau^0_k-\varepsilon,\tau^0_k+\varepsilon)$ the mapping $z(t,t_1,z^1,\mu)$ is $C^2$~-- smooth with respect to the variables $(t_1,z^1,\mu)\in U$. Moreover, this solutions have exactly $N$ impact instants
$\tau_j(t_1,z^1,\mu_1)$ $(j=1,\ldots,N)$ over the segment $[t_-,t_+]$. These instants and corresponding velocities $Y_j=-y_1(\tau_j(t_1,z^1,\mu_1)-0,t_1,z^1,\mu_1)$ $C^2$ smoothly depend on $t_1,z^1$ and $\mu_1$.
\end{lemm}

\textbf{Proof.} Let the number $k$ be such that $t_0\in [\tau_k,\tau_{k+1}]$ (assume, if necessary $\tau_0=t_-$, $\tau_{N+1}=t_+$). The solution $z(t)$ of the vibroimpact system \1 is also one of \eqref{mgr1.1} over the segment $[\tau_k,\tau_{k+1}]$. The impact instants $\tau_k$ and $\tau_{k+1}$ as well as the impact points $$z(\tau_k-0,t_1,z^1,\mu_1)\qquad \mbox{and} \qquad z(\tau_{k+1}-0,t_1,z^1,\mu_1)$$ smoothly depend on their parameters. Similarly, the values $\tau_{k-1}$ and $z(\tau_{k-1}-0,t_1,z^1,\mu_1)$ smoothly depend on $\tau_k$ and $z(\tau_k-0,t_1,z^1,\mu_1)$, as well as $\tau_{k+2}$ and $z(\tau_{k+2}-0,t_1,z^1,\mu_1)$ are $C^1$~-- smooth functions of $\tau_{k+1}$ and $z(\tau_{k+1}-0,t_1,z^1,\mu_1)$ and so on. This proves the Lemma \ref{leg1}.

\begin{cond}\label{cog2}(fig.\,1.). There exists a continuous family of $T$~-- periodic solutions $$\varphi(t,\mu)=\col(\varphi^x_1(t,\mu),\varphi^y_1(t,\mu),\ldots,\varphi_n^x(t,\mu),\varphi_n^y(t,\mu)) \qquad t\in{\mathbb R},\quad \mu\in J$$ of the system \1, satisfying the following properties.
\begin{enumerate}
\item For all $\mu\geqslant 0$ the component $\varphi_1^x(t,\mu)$ has exactly $N+1$ zeros $\tau_0(\mu),\ldots,\tau_{N}(\mu)$ over the period $[0,T)$.
\item The velocities $Y_k(\mu)=-\varphi_1^y(\tau_k(\mu)-0,\mu))$ are such that
\begin{equation}\label{mgr2.1} \begin{array}{l}
Y_0(\mu)>0 \qquad\mbox{for all}\quad \mu>0,\qquad Y_0(0)=0,\\
\qquad f_1(\tau_0(0),0,0,\bar\varphi(\tau_0(0),0),0)=\phi_0>0,\\
Y_k(\mu)>0, \qquad \forall \mu\in J, \quad k=1,\ldots, N.
\end{array} \end{equation}
Here $\bar\varphi(t,\mu)=\col(\varphi^x_2(t,\mu),\varphi^y_2(t,\mu),\ldots,\varphi_n^x(t,\mu),\varphi_n^y(t,\mu))$.
\vspace{30pt}
\begin{figure}[!ht]
\begin{center}
\epsfxsize=1.22in \epsfysize=1.86in \epsfbox{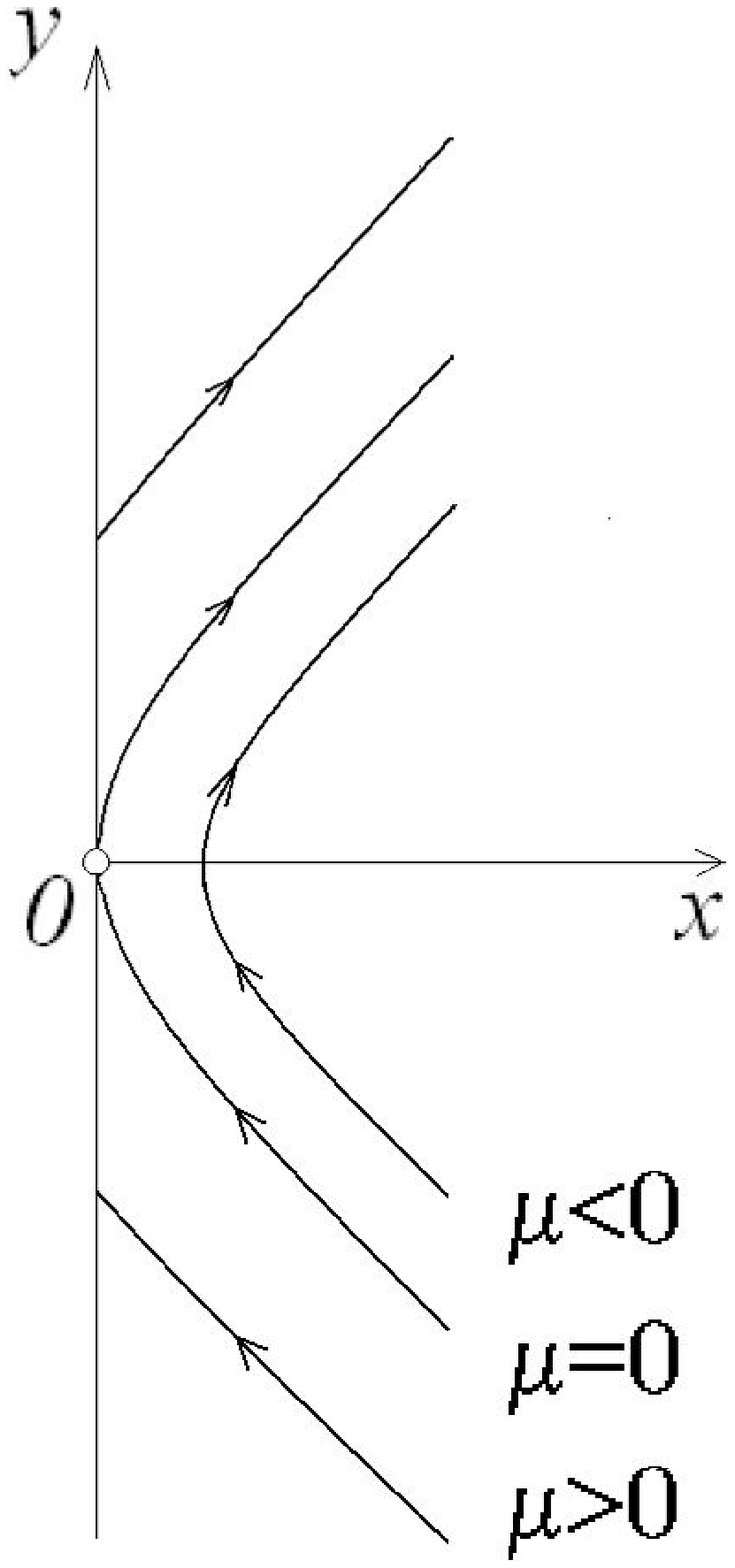}
\vspace{-5mm} {\small \noindent {\it Fig.\, 1.}}
\end{center} \end{figure}
\item The instants $\tau_k(\mu)$ and the velocities $Y_k(\mu)$ continuously depend on $\mu\in J$.
\end{enumerate}\end{cond}

Denote
$$\begin{array}{c}
\varphi^x(t,\mu)=\col(\varphi^x_1(t,\mu),\ldots,\varphi_n^x(t,\mu)), \qquad \varphi^y(t,\mu)=\col(\varphi^y_1(t,\mu),\ldots,\varphi_n^y(t,\mu)),\\
\bar{\varphi}^x(t,\mu)=\col(\varphi^x_2(t,\mu),\ldots,\varphi_n^x(t,\mu)), \qquad \bar{\varphi}^y(t,\mu)=\col(\varphi^y_2(t,\mu),\ldots,\varphi_n^y(t,\mu)).\end{array}$$
We may suppose without loss of generality that $\tau_0(\mu)\equiv 0$, $\bar\varphi (0,0)=0$. This may be obtained by the transformation $t'=t-\tau_0(\mu)$. Define
$$\theta_0=\dfrac12 \min\left(\min_{\mu\in J} \tau_1(\mu),T-\max_{\mu\in J} \tau_N(\mu)\right).$$ Fix a small positive $\theta<\theta_0$ and consider the shift mapping for the system \1, given by the formula $S_{\mu,\theta}(z^0)=z(T-\theta+ 0,-\theta,z^0,\mu)$. For small positive $\mu$ and $\theta$ the mapping $S_{\mu,\theta}$ is $C^1$~- smooth in a neighborhood of the point $z_{\mu,\theta}=\varphi(-\theta,\mu)$. Without loss of generality, we may assume that $$\lim\limits_{\mu,\theta\to 0+} z_{\mu, \theta}=0.$$

\noindent\textbf{3. The separatrix.} Denote
$\Gamma_{\mu,\theta}=\{z^0\in \Lambda:\exists t_1\in [-T,T]: z_1(t_1,-\theta,z^0,\mu)=0\}.$

\begin{lemm}\label{leg3} There exists a neighborhood $U_0$  of zero such that if the parameters $\mu$ and $\theta$ are small enough, the set $\Gamma_{\mu,\theta}\bigcap U_0$ is a surface of the dimension $2n-1$, which is the graph of the $C^2$ smooth function $x_1=\gamma_{\mu,\theta}(\bar x,y)$ (fig.\,2). Moreover,
\begin{equation}\label{mgr3.2}
\gamma_{\mu,\theta}(z)=y^2_1\left(\dfrac1{f_1(-\theta,z_{\mu,\theta},\mu)}+\widetilde{\gamma_{\mu,\theta}}(z)\right),
\end{equation}
where $\widetilde{\gamma_{\mu,\theta}}$ is a $C^2$ smooth function such that $\widetilde{\gamma_{\mu,\theta}}(0)=0$.
\vspace{30pt}
\begin{figure}[!ht]
\begin{center}
\epsfxsize=1.83in \epsfysize=2.32in \epsfbox{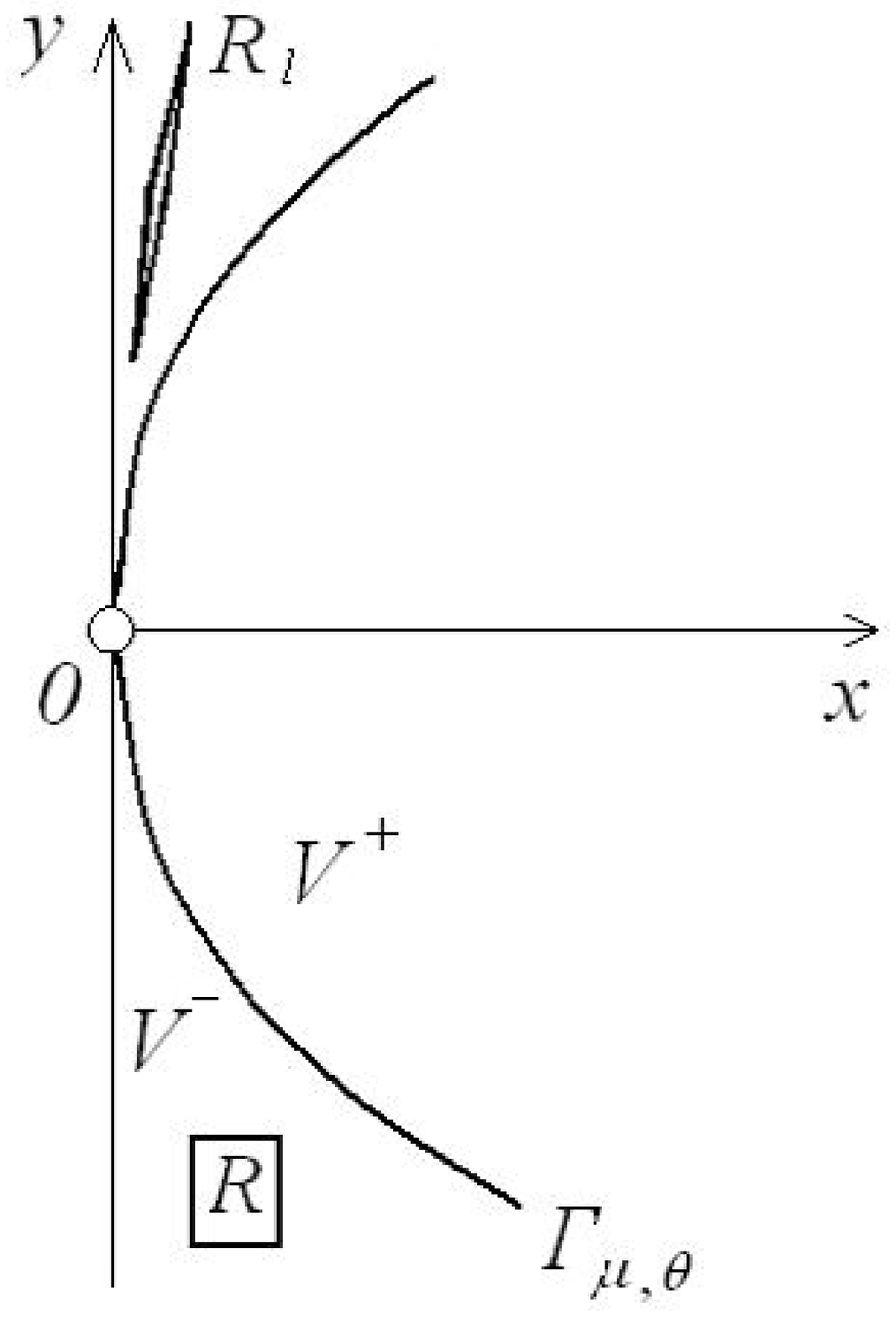}
\vspace{-5mm} {\small \noindent {\it Fig.\, 2.}}
\end{center} \end{figure}\end{lemm}

\textbf{Proof}. Take a point $\zeta\in \Gamma_{\mu,\theta}$. Let the instant $t_0$ be such that $z_1(t_0,-\theta,\zeta)=0$, $s=t-t_0$, $z_1(t+0,-\theta,\zeta)=\col(x_1(t),y_1(t))$. Let us show that if $t_0$ is close enough to $-\theta$, we may take $s_0\geqslant |t_0+\theta|$ so that the function $x_1(t_0+s)$ does not have zeros on $[-s_0,s_0]$, except $s=0$. Otherwise, there exists a sequence $t^k_0\to -\theta$ (suppose without loss of generality, that $t_0^k>-\theta$ and the sequence decreases), a sequence $t_1^k\in [-\theta, t_0^k)$ and one, consisting of solutions, uniformly bounded on the segment $[-\theta,\max t_0^k]$:
$$z^k(t)=\col (z^k_1(t),\ldots, z^k_n(t))=\col(x^k_1(t),y^k_1(t),\ldots,x^k_n(t),y^k_n(t))$$€
of the system \1, such that $z^k_1(t_0^k)=0$, $x^k_1(t_1^k)=0$.
\vspace{30pt}
\begin{figure}[!ht]
\begin{center}
\epsfxsize=1.6in \epsfysize=1in \epsfbox{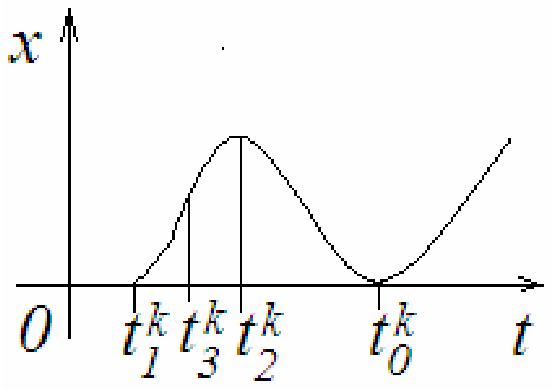}
\vspace{-5mm} {\small \noindent {\it Fig.\, 3.}}
\end{center} \end{figure}
Also, there exist time instants $t_2^k\in (t_1^k,t_0^k)$, that
$\dot x^k_1(t_2^k)=0$ and instants $t_3^k\in (t_2^k,t_0^k)$ such that $\ddot x^k_1(t_3^k)=0$. Moreover, $t_3^k\to -\theta$, $x^k_1(t_3^k)\to 0$, $\dot x^k_1(t_3^k)\to 0$. Without loss of generality, we assume that $\bar z(t^k_3)\to \bar z^0$. Then $\ddot x_1(t_3^k)\to f_1(-\theta,{\bar z}^0,\mu)=0$. This contradicts to \eqref{mgr2.1}.

Then for all $s\in [-s_0,s_0]$ the function $x_1(t_0+s)$ can be presented as series
\begin{equation}\label{mgr3.3}
  x_1(t_0+s)=X_2s^2+X_3s^3+\ldots
\end{equation}
Differentiating \eqref{mgr3.3}, we obtain that $\dot x_1(t_0+s)=2X_2s+3X_3s^2+\ldots$. On the other hand,
$X_2=\ddot x_1(t_0+0)/2\to f_1(-\theta,z_{\mu,\theta},\mu)/2$ as $t_0\to -\theta$. Then $$x_1(-\theta)=f_1(-\theta,z_{\mu,\theta},\mu)(t_0+\theta)^2(1+o(1))/2;\qquad y_1(-\theta)=f_1(-\theta,z_{\mu,\theta},\mu)(t_0+\theta)(1+o(1)).$$ Since $y_1=\dot x_1$, the formula \eqref{mgr3.2} is true. The lemma is proved.

Take a small parameter $\varsigma>0$ such that the sets
$V_{\mu,\theta}=\{z\in \Lambda: \|z-z_{\mu,\theta}\|\leqslant \varsigma\}\subset U_0$,
  $$V_{\mu,\theta}^-=\{(x,y)\in V_{\mu,\theta}: x_1<\gamma_{\mu,\theta}(\bar x,y)\}, V_{\mu,\theta}^+=\{(x,y)\in V: x_1>\gamma_{\mu,\theta}(\bar x,y)\}$$
are correctly defined and nonempty.

\noindent\textbf{4. The main result.} Consider the matrix
$$A=\lim\limits_{\theta\to 0+}\left.\pd{z}{z^0}(T-\theta+0, \theta,z^0,0)\right|_{z^0=z_{0,\theta}}.$$
Let $\Delta_0=\det A$. Denote the elements of the matrix $A$ by $a_{ij}$ and ones of the matrix $A^{-1}$ by $\alpha_{ij}$. Denote the columns of matrices $A$ and $A^2$ by $A_j$ and $A_j^2$ respectively, the strings of the matrix $A^{-1}$ by ${\cal A}_j$. If $n>1$ and $a_{12}\neq 0$ consider the $(2n-2)\times (2n-2)$ matrix
$\bar A=({\bar a}_{ij})$, defined by formulae $${\bar a}_{ij}=\dfrac{-a_{1j+2}a_{i+2\,2}}{a_{12}}+a_{i+2\,j+2},\qquad (i,j=1,\ldots, 2n-2).$$ Similarly, if $n=1$ and $\alpha_{12}\neq 0$, we define the matrix $\overline{\cal A}=(\bar\alpha_{ij})$
$${\bar \alpha}_{ij}=\dfrac{-\alpha_{1j+2}\alpha_{i+2\,2}}{\alpha_{12}}+\alpha_{i+2\,j+2},\qquad (i,j=1,\ldots, 2n-2).$$

Assume that the at least one of the following statements is true.

\begin{cond}\label{cog3} \
\begin{enumerate}
\item Either $n=1$ or the matrix $\bar A$ does not have eigenvalues on the unit circle in $\aC$.
\item \begin{equation}\label{mgr4.1} a_{12}>0,\qquad \sum_{k=1}^{2n} a_{1k}a_{k2}<0\end{equation}
\end{enumerate}
\end{cond}

\begin{cond}\label{cog4}\
\begin{enumerate}
\item Either $n=1$ or the matrix $\bar {\cal A}$ does not have eigenvalues on the unit circle in $\aC$.
\item \begin{equation}\label{mgr4.2} \alpha_{12}>0,\qquad \sum_{k=1}^{2n} \alpha_{1k}\alpha_{k2}<0. \end{equation}
\end{enumerate}
\end{cond}

From the geometrical point of view, the first items of Conditions 3 and 4 mean that the fixed points $z_{\mu,\theta}$ are saddle hyperbolic and the inequalities \eqref{mgr4.1} and \eqref{mgr4.2} provide that the corresponding stable and unstable manifolds intersect. This will be shown below.

Later on we shall suppose that Condition \ref{cog3} is satisfied. Otherwise, we consider the mapping $S^{-1}_{\mu,\theta}$ instead of $S_{\mu,\theta}$. Then the matrix $A$ is replaced with $A^{-1}$, and the condition \eqref{mgr4.2} with \eqref{mgr4.1}. The similar reasonings shall prove the statement of the lemma in the considered case (see the right part of the figure 4). All the proofs given below, may be repeated for this case.

\vspace{30pt}
\begin{figure}[!ht]
\begin{center}
\epsfxsize=3.66in \epsfysize=2.55in \epsfbox{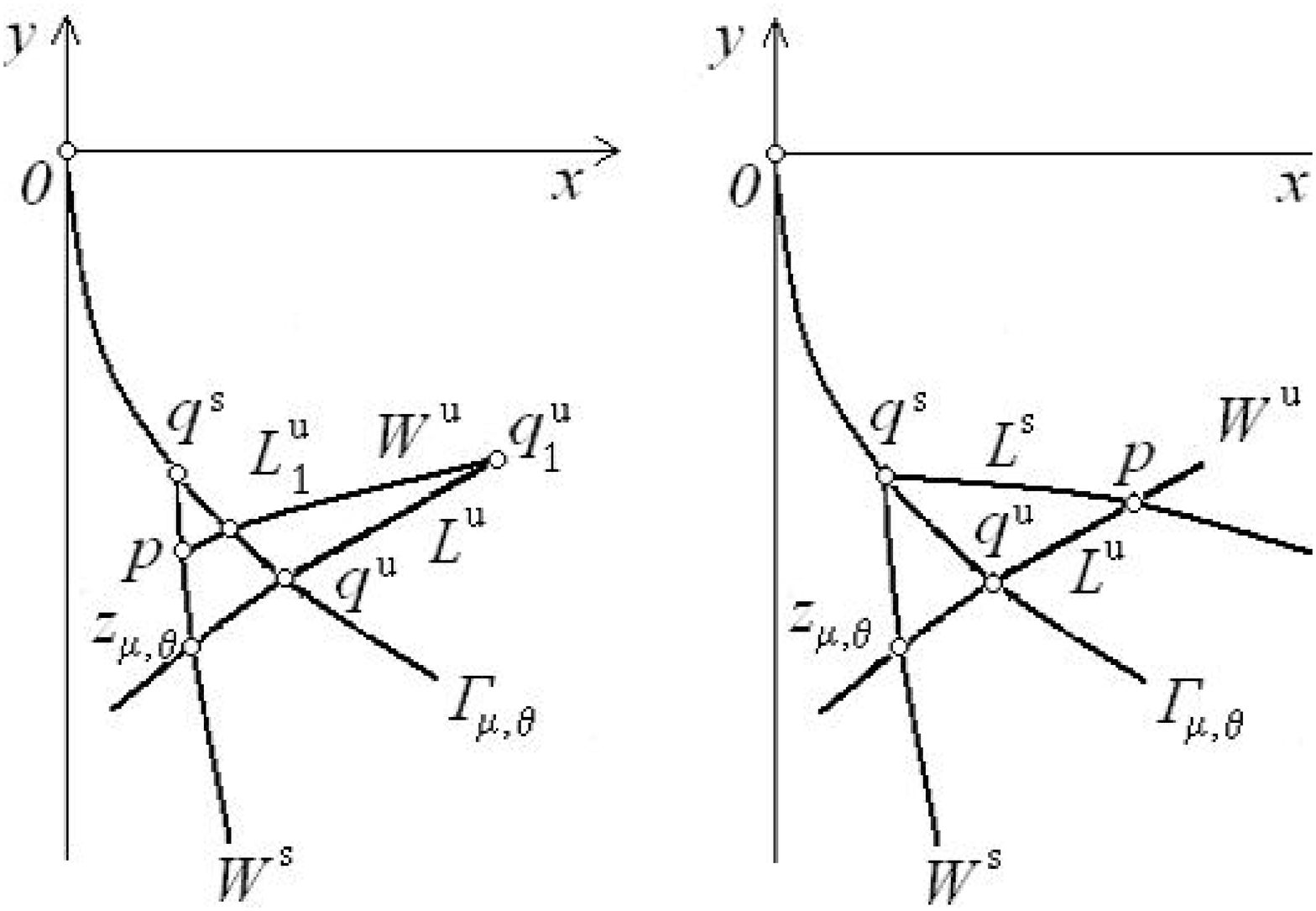}
\vspace{-5mm} {\small \noindent {\it Fig.\, 4.}}
\end{center} \end{figure}

\begin{theo}\label{thg1} If Condition \ref{cog2} and one of Conditions \ref{cog3} or \ref{cog4} are satisfied, there exist values $\mu_0>0$ and $\theta>0$ such that for all $\mu\in (0,\mu_0)$ the mapping $S_{\mu,\theta}$ is chaotic in the sense of \cite{32.5.}. More precisely, there exists an integer $m$ and a compact set $K=K_{\mu,\theta}$ invariant with respect to $S^m_{\mu,\theta}$ such that the following conditions are satisfied.
\begin{itemize}
\item[\bf 1)] There exists a neighborhood $U_{\mu,\theta}$ of the set $K_{\mu,\theta}$ such that the mapping
$S^m_{\mu,\theta}|_{U_{\mu,\theta}}$ is a diffeomorphism. The invariant set $K_{\mu,\theta}$ is hyperbolic.
\item[\bf 2)] The mapping $S^m_{\mu,\theta}|_{K_{\mu,\theta}}$ has infinitely many periodic points.
\item[\bf 3)] The periodic points of $S^m_{\mu,\theta}$ are dense in $K_{\mu,\theta}$.
\item[\bf 4)] The set $K_{\mu,\theta}$ is transitive i.e. there exists a point $p_{\mu,\theta}\in K_{\mu,\theta}$, such that the orbit $\{S_{\mu,\theta}^{km}(p_{\mu,\theta})\}_{k=-\infty}^\infty$ is dense in $K_{\mu,\theta}$.
\end{itemize}
\end{theo}

\begin{rema} The similar results may be obtained for the systems with any finite number of grazings over the period. \end{rema}

\noindent\textbf{5. Grazing.} Now we start to prove the theorem. Note that all the mappings $S_{\mu,\theta}$, corresponding to the same value of $\mu$, are conjugated. Fix a small value $\mu>0$ and a solution
$$z^0(t)=\col(x^0_1(t),y^0_1(t),\ldots,x^0_n(t),y^0_n(t))$$
of the corresponding system, having an impact at the instant $t_0$. Suppose that the corresponding normal velocity $Y_{01}=-y_1^0(t_0-0)$ is nonzero. We consider $Y_{01}$ as a small parameter. Denote $\bar Z_0=\bar z^0(t_0-0)$. Fix a positive value $s_0$ and consider the mapping $$G(\zeta)=z(t_0+s_0,t_0-s_0,\zeta,\mu),$$ defined in a neighborhood of the point $\zeta^0=z^0(t_0-s_0)$. Here we assume that the point $\zeta^0$ and the parameter $s_0$ are chosen so that there exists a neighborhood $\Omega\ni \zeta^0$ such that any solution $z(t)=z(t,t_0-s_0,z_-,\mu)$ ($z_-\in \Omega$) has exactly one impact over the segment $[t_0-s_0,t_0+s_0]$. Denote the corresponding instant by $t_1=t_1(z_-)$ and the normal velocity, defined similarly to $Y_{01}$, by $Y_1=Y_1(z_-)$. Let $\bar Z=\bar z(t_1)$ be the tangent component of the solution $z(t)$ at the impact instant. Take the values $s_\pm=s_\pm(z_-)$ so that $t_0\pm s_0=t_1(z_-)\pm s_\pm(z_-)$ for all $z_-\in \Omega$. The mapping $G$ is smooth in the neighborhood of the point $\zeta_0$, let us estimate the Jacobi matrix $DG(\zeta_0)$. Denote
$$\begin{array}{l} z_+=z(t_0+s_0)=z(t_0+s_0,t_1+0,0,rY_1,\bar Z,\mu),\\
x_+=x(t_0+s_0)=x(t_0+s_0,t_1+0,0,rY_1,\bar Z,\mu),\\
y_+=z(t_0+s_0)=y(t_0+s_0,t_1+0,0,rY_1,\bar Z,\mu),\\
x_-=x(t_0-s_0)=x(t_0-s_0,t_1-0,0,-Y_1,{\bar Z},\mu),\\
y_-=y(t_0-s_0)=y(t_0-s_0,t_1-0,0,-Y_1,{\bar Z},\mu).\end{array}$$

Similarly, we define the values $x_{1,\pm}$, $y_{1,\pm}$, ${\bar x}_{\pm}$, ${\bar y}_{\pm}$. Consider the Taylor formula for values $z_\pm$ as functions of $s_\pm$:
\begin{equation}\label{mgr5.0}\begin{array}{l}
x_{1,-}=Y_1s_-+f_1(t_1,0,-Y_1,{\bar Z},\mu)s_-^2/2+\rho_{x1-}(s_-,t_1,Y_1,{\bar Z},\mu)s_-^3;\\
y_{1,-}=-Y_1-f_1(t_1,0,-Y_1,{\bar Z},\mu)s_-+\rho_{y1-}(s_-,t_1,Y_1,{\bar Z},\mu)s_-^2;\\
{\bar x}_-={\bar x}(t_1)-{\bar y}(t_1-0)s_-+{\bar f}(t_1,0,-Y_1,{\bar Z},\mu)s_-^2/2+\rho_{{\bar x}-}(s_-,t_1,Y_1,{\bar Z},\mu)s_-^3;\\
{\bar y}_-={\bar y}(t_1-0)-{\bar f}(t_1,0,-Y_1,{\bar Z},\mu)s_-+\rho_{{\bar y}-}(s_-,t_1,Y_1,{\bar Z},\mu)s_-^2;\\
x_{1,+}=rY_1s_++f_1(t_1,0,rY_1,{\bar Z},\mu)s_+^2/2+\rho_{x1+}(s_+,t_1,Y_1,{\bar Z},\mu)s_+^3;\\
y_{1,+}=rY_1+f_1(t_1,0,rY_1,{\bar Z},\mu)s_++\rho_{y1+}(s_+,t_1,Y_1,{\bar Z},\mu)s_+^2;\\
{\bar x}_+={\bar x}(t_1)+{\bar y}(t_1-0)s_++{\bar f}(t_1,0,rY_1,{\bar Z},\mu)s_+^2/2+\rho_{{\bar x}+}(s_+,t_1,Y_1,{\bar Z},\mu)s_+^3;\\
{\bar y}_+={\bar y}(t_1-0)+{\bar f}(t_1,0,rY_1,{\bar Z},\mu)s_++\rho_{{\bar y}+}(s_+,t_1,Y_1,{\bar Z},\mu)s_+^2.\\
\end{array}\end{equation}
Here all functions, denoted by the letter $\rho$ with different indices, are $C^2$ smooth with respect to all arguments except $s_\pm$. Denote
$$f_{0k+}=f_k(t_0,0,rY_{01},{\bar Z}_0,\mu), \qquad f_{0k-}=f_k(t_0,0,-Y_{01},{\bar Z}_0,\mu), \qquad
{\widetilde r}=r+\left.\dfrac{\partial r}{\partial y_1}\right|_{y_1=Y_{01}}Y_{01}.$$
It follows from \eqref{mgr5.0} that
$$\begin{array}{rl}
  &\left.\pd{z_+}{(s_+,Y_1,\bar Z)}\right|_{s_+=0,\, Y_1=Y_{01},\, \bar Z=\bar Z_0}=
  \begin{pmatrix}
    rY_{01} & 0 & 0\\
    f_{01+} & \widetilde{r} & 0\\
    Q_+& 0 & E_{2n-2}
  \end{pmatrix};  \\[5pt]
  &\left.\pd{z_-}{(s_-,Y_1,\bar Z)}\right|_{s_-=0,\, Y_1=Y_{01},\, \bar Z=\bar Z_0}=
  \begin{pmatrix}
    Y_{01} & 0 & 0\\
   -f_{01-} & -1 & 0\\
   Q_- & 0 & E_{2n-2}
  \end{pmatrix}.
\end{array}$$
Here $E_{2n-2}$ is the unit matrix of the corresponding size,
$$\begin{array}{c} Q_+=\col(y_2(t_0-0), f_{02+}, \ldots,y_n(t_0-0), f_{0n+}),\\
 Q_-=\col(-y_2(t_0-0), -f_{02+}, \ldots,-y_n(t_0-0), -f_{0n+}).\end{array}$$

Denote
$$f'_k=\left.\dfrac{\partial f_k(t_0,0,y_1,\bar Z_0,\mu)}{\partial y_1}\right|_{y_1=0}$$
Clearly, $ds_+/ds_-=-1$. Then, similarly to the results of the paper \cite{48.}, we obtain
\begin{equation}\label{mgr5.2}
B=\lim_{s_\pm \to 0}\pd{z_+}{z_-}=\begin{pmatrix}
-r & 0 & 0 \\
b_{21} & -\widetilde{r} & 0\\
\bar{B}_1& 0& E_{2n-2}
\end{pmatrix}\end{equation}
(fig.\,2). Here $\bar{B}_1=\col(b_{31}, \ldots b_{2n\,1}),$
$$\begin{array}{c}
b_{21}=-(f_{01+}+\widetilde{r}f_{01-})/Y_{01}=-(r+1)\phi_0(1+O(Y_{01}))/Y_{01},\quad b_{2j-1\,1}=0,\\
b_{2j\,1}=(f_{0j+}-f_{0j-})/Y_{01}=(r+1)f'_k+O(Y_{01}).\end{array}$$
Note that, $\det B=r\widetilde{r}\to r^2$ as $Y_{01}\to 0$.

It is shown on the figure 2 how the small neighborhood $R$ of the form
$$R=\{(\xi_1,\psi_1,\ldots,\xi_n,\psi_n):|\xi_j-\xi^0_j|\leqslant \delta,\,|\psi_j-\psi^0_j|\leqslant \delta\}$$
is stretching under the action of the mapping $S_{1,\mu,\theta}$, defined by the formula $$S_{1,\mu,\theta}(\zeta)=z(\theta+ 0,-\theta,\zeta,\mu).$$ Here $R_1=S_{1,\mu,\theta}(R)$.

\noindent\textbf{6. Lyapunov exponents.} We check that the fixed point $z_{\mu,\theta}$ of the mapping $S_{\mu,\theta}$ is hyperbolic and estimate the bigger and the smaller absolute value of eigenvalues of the matrix $D=DS_{\mu,\theta}(z_{\mu,\theta})$ and ones of small perturbations of this matrix. For a fixed value $\mu>0$ denote $D=DS_{\mu,\theta}(z_{\mu,\theta})$. The mapping $S_{\mu,\theta}$ can be presented as the composition $S_{\mu,\theta}=S_{2,\mu,\theta}\circ S_{1,\mu,\theta}$, where $S_{2,\mu,\theta}(\zeta)=z(T-\theta+0,\theta,\zeta,\mu)$. The Jacobi matrix $A_{\mu,\theta}=DS_{2,\mu,\theta}(S_{1,\mu,\theta}(z_{\mu,\theta}))$ tends to $A$ as
$\mu,\theta\to 0$. The matrix $B_{\mu,\theta}=DS_{1,\mu,\theta}(z_{\mu,\theta})$ is of the form \eqref{mgr5.2}, where $Y_{01}=Y_0(\mu)$. Then
$$D=A_{\mu,\theta}B_{\mu,\theta}=(-(r+1)A_2\phi_0(1+O(Y_0))/Y_0,-rA_2(1+O(Y_0)),A_3+O(Y_0),\ldots,A_{2n}+O(Y_0)).$$
Since $\det D=\det A_{\mu,\theta} \det B_{\mu,\theta}=(r^2+O(Y_0))\Delta_0$, if the inequalities \eqref{mgr4.1} take place and if $\mu$ is small enough, one of the eigenvalues of the matrix $D$ is $$\lambda_+=-(r+1)a_{12}\phi_0(1+O(Y_0))/Y_0.$$ The corresponding eigenvector $u_+$ equals to $A_2+O(Y_0)$. The eigenvalue $\lambda_+$ is of the multiplicity 1, the linear space, corresponding to other eigenvalues, tends as $\mu\to 0$ to the hyperplane $\pi_1$, given by the condition $x_1=0$. Since $a_{12}\neq 0$, the vector $A_2$ is out of $\pi_1$. The matrix $D^{-1}$ satisfies the following asymptotic estimate
$$D^{-1}=\dfrac1{r^2}\begin{pmatrix}
-r{\cal A}_1+O(Y_0)\\
(r+1){\cal A}_1(1+O(Y_0))/Y_0\\
{\cal A}_3+O(Y_0)\\
\dots\\
{\cal A}_{2n}+O(Y_0)
\end{pmatrix}.$$
It follows from the form of this matrix, that the matrix $D^{-1}$ has the eigenvalue $$\lambda_-^{-1}=(r+1)\phi_0 \alpha_{12}(1+O(Y_0))/(r^2Y_0).$$ The corresponding eigenvector satisfies the asymptotical estimate $u_-=e_{2}+O(Y_0)$.

If $n=1$, it is clear that the matrix $AB$ as well as matrices $DS_{\mu,\theta}(z)$ corresponding to points $z$ of a small neighborhood of $z_{\mu,\theta}$ are hyperbolic. Otherwise we need the following lemma.

\begin{lemm}\label{leev} If Condition 3 is satisfied there exist positive constants $\mu_0$ and $\varepsilon_0$ such that for any matrix $A'$ such that $\|A'-A\|<\varepsilon_0$ and any $\mu\in (0,\mu_0)$ the matrix $A'B_{\mu,\theta}$ does not have eigenvalues on the unit circle in $\aC$.
\end{lemm}

\textbf{Proof}. Suppose the statement of the lemma is not true. Then there exist a sequence of matrices $A_k\to A$ and sequences $\mu_k,\theta_k\to 0$ such that all matrices $A_kB_{\mu_k,\theta_k}$ have eigenvalues of the form $\lambda_k=\exp(i\psi_k)$. Denote the corresponding eigenvectors by
$u^k=\col(u^k_1,\ldots, u^k_{2n})$. Without loss of generality we may assume that $\|u^k\|=1$ for all $k$ and that
$$\lambda_k\to \lambda_0=\exp(i\psi_0),\qquad u^k\to u^0=\col(u^0_1,\ldots, u^0_{2n})$$ as $k\to\infty$. Denote the elements of the matrices $A_k$ and $B_{\mu_k,\theta_k}$ by $a^k_{ij}$ and $b_{ij}^k$ respectively ($k\in \aN$, $i,j=1,\ldots, 2n$). Let $r_k=r(Y_0(\mu_k),\mu_k)$, $\widetilde{r}_k=\widetilde{r}(Y_0(\mu_k),\mu_k)$,
$$d^k_i=\sum_{j=1}^{2n}a^k_{ij}b^k_{j1}.$$
Note that since $a_{12}^k>a_{12}/2>0$, $b_{12}^k\to \infty$ and other values $b_{12}^k$ are uniformly bounded, the sequence $d_1^k$ tends to $+\infty$ as $k\to+\infty$ and
$$\lim_{k\to+\infty}d_i^k/d_1^k=a_{i2}/a_{12}, \qquad (i=1,\ldots,2n).$$ From the definition of eigenvalues and eigenvectors we obtain
\begin{equation}\label{eqev1}
d^k_ju_1^k-\widetilde{r}_ka^k_{j2}u^k_2+\sum_{s=3}^{2n}a^k_{js} u_s^k=\lambda_k u_j^k \qquad j=1,\ldots, 2n.
\end{equation}

It follows from the first equation \eqref{eqev1} that
$$u_1^k=\dfrac{-\widetilde{r}_ka_{12}u^k_2+\sum\limits_{s=3}^{2n} a_{1s}^k u^k_s}{-d^k_1+\lambda_k}\to 0$$
as $k\to +\infty$. Consequently, $u_1^0=0$. Substituting this expression to Eq. \eqref{eqev1}, corresponding to $j>1$, we have
$$\dfrac{d^k_j(-\widetilde{r}_ka^k_{12}u^k_2+\sum\limits_{s=3}^{2n} a_{1s}^k u^k_s)}{-d^k_1+\lambda_k}-\widetilde{r}_ka^k_{j2}u^k_2+\sum_{s=3}^{2n}a^k_{js} u_s^k=\lambda_k u_j^k \qquad j=2,\ldots, 2n.$$

Proceeding to the limit as $k\to\infty$ we obtain
$$\sum_{s=3}^{2n} \left(\dfrac{-a_{1s}a_{j2}}{a_{12}}+a_{js}\right) u^0_s= \lambda_0u_j^0, \qquad j=2,\ldots,2n.$$
Clearly, at least one of values $u_j^0$ ($j=3,\ldots,2n$) is nonzero. Then $\lambda_0$ is an eigenvalue of the matrix $\bar A$. This contradicts to our assumptions.

\noindent\textbf{7. Homoclinic point.} The mapping $S_{\mu,\theta}$ is differentiable at the points of the set $V^-_{\mu,\theta}$. Due to Lemma \ref{cog3} the eigenvalues of Jacobi matrices $D$ are out of the unit circle, provided $\mu$ is small. Then, due to the Perron theorem, in a small neighborhood of the point $z_{\mu,\theta}$ there exist the local stable manifold $W^s$ and the unstable one $W^u$ of the mapping $S_{\mu,\theta}$. Both of them are smooth surfaces. Let $M^{s,u}=T_{z_{\mu,\theta}} W^{s,u}$, $n_s=\dim W^s$, $n_u=\dim W^u=2n-n_s$. Select the orthonormal basises $e^s_1, \ldots, e^s_{n_s}$ and $e^u_1, \ldots, e^u_{n_u}$ in the spaces $M^u$ and $M^s$ respectively so that $e^s_1=u_-/\|u_-\|$ and $e^u_1=u_+/\|u_+\|$. Extend the stable and unstable manifolds up to the invariant sets $W^s$ and $W^u$ of the mapping $S_{\mu,\theta}$. The obtained sets consist, generally speaking, of a countable number of the connected components. Every of these components is a partially smooth manifold.

\begin{lemm}\label{leg4} The manifolds $W^s$ and $W^u$ intersect transversally at a point $p\neq z_{\mu,\theta}$ (fig.\, 4). \end{lemm}

\textbf{Proof.} If conditions \eqref{mgr4.1} be satisfied then, for small values of the parameter $\mu$ the manifold $W^u$ intersects transversally the surface $\Gamma_{\mu,\theta}\bigcap V_{\mu,\theta}$. Denote the manifold, obtained in the intersection, by $q^u$. The dimension of $q^u$ equals to $n_u-1$. Denote $q^u_1=S_{\mu,\theta}(q^u)$. The neighborhood $V_{\mu,\theta}$ of the point $z_{\mu,\theta}$ may be chosen so that both the manifolds $q^u$ and $q^u_1$ intersect with $V_{\mu,\theta}$ and the intersection $W^u\bigcap V_{\mu,\theta} \setminus q^u$ consists of two connected components. Denote one, which does not contain the point $z_{\mu,\theta}$, by $L^u$. The neighborhood $V_{\mu,\theta}$ may be chosen so that $\diam V_{\mu,\theta}\to 0$ as $\mu,\theta\to 0+$. Let $L^u_1=S_{\mu,\theta}(L^u)$. For any $z\in L^u$ the tangent space $M^u(z)=T_z W^u$ is the linear hull of unit orthogonal vectors $E^u_1(z),E^u_2(z), \ldots, E^u_{n_u}(z)$, which can be chosen so that for any $k=1,\ldots,n_u$ $$\lim_{\mu,\theta\to 0+}\max_{z\in L^u}\|E^u_k(z)-e^u_k\|=0.$$
The surface $W^u$ is not smooth in the neighborhood of the manifold $q^u_1$. For the points $z\in L^u_1$ the tangent space $M^u(z)=T_z W^u$ is the linear hull of unit vectors $\widetilde{E}^u_1(z),E^u_2(z), \ldots, E^u_{n_u}(z)$, where
$$\lim_{\mu,\theta\to 0+}\max_{z\in L^u_1}\|E^u_k(z)-e^u_k\|=0$$
for all $k>1$ and
$$\lim_{\mu,\theta\to 0+}\max_{z\in L^u_1}\|\widetilde{E}^u_1(z)-A^2_2/\|A^2_2\|\|=0.$$
It follows from the conditions \eqref{mgr4.1}, that for any $z\in L^u_1$ the set $\widetilde{E}^u_1(z),E^u_2(z), \ldots, E^u_{n_u}(z)$ is linearly independent. The vectors $A_2$ and $A_2^2$ lie at different half-spaces, separated by the hyperplane $\pi_1$. Consequently, for all
$z_0\in L^u$, $z_1\in L^u_1$ and $z_2\in W^u\bigcap V_{\mu,\theta}$ the vectors $E^u_1(z_0)$ and $\widetilde{E}^u_1(z_1)$ lie in different half-spaces, separated by the linear hull of vectors $E^u_2(z_2), \ldots, E^u_{n_u}(z_2),E^s_1(z_2), \ldots, E^s_{n_s}(z_2)$. This means that for the small values of the parameter $\mu$ and $\theta$ the neighborhood $V_{\mu,\theta}$ may be chosen so that the surfaces $L^u_1$ and $W^s$ intersect transversally. This proves the lemma for the considered case (see the left part of the figure 4).

The Smale-Birkhof theorem \cite{76.} on the existence of a chaotic invariant set in a neighborhood of a homoclinic point is not applicable in the considered case since the mapping $S_{\mu,\theta}$ is discontinuous. However, the similar techniques will help us to find a chaotic set of the mapping $S_{\mu,\theta}$.

\textbf{Acknowledgements.} This work was supported by the UK Royal Society, by the Russian Federal Program "Scientific and pedagogical cadres", grant no. 2010-1.1-111-128-033 and by the Chebyshev Laboratory (Department of Mathematics and Mechanics, Saint-Petersburg State University) under the grant \\ 11.G34.31.2006 of the Government of the Russian Federation.

\newpage \thispagestyle{empty}
ABSTRACT

{Kryzhevich S.\,G.}

{\bf Smale horseshoes and grazing bifurcation in vibroimpact systems}

Bifurcations of dynamical systems, described by a second order differential equations and by an impact condition are studied. It is shown that the variation of parameters when the number of impacts of a periodic solution increases, leads to the occurrence of a hyperbolic chaotic invariant set.

\end{document}